\documentclass[10pt,a4paper]{article}
\usepackage{amsmath,amsthm,amsfonts,amssymb}
\usepackage{url}

\usepackage[usenames,dvipsnames]{color}
\newenvironment{red}{\color{red}}{}
\newtheorem{theorem}{Theorem}

\newtheorem{pro-example}[theorem]{Example}

\newtheorem{pro-remark}[theorem]{Remark}
\newenvironment{remark}{\begin{pro-remark}\rm}{\cqfd\end{pro-remark}}

\newtheorem{pro-definition}[theorem]{Definition}

\def\Schreier{\textsf{Schreier}}

\def\calP{\mathcal{P}}

\let\phi\varphi

\let\epsilon\varepsilon

\newcommand{\inv}{^{-1}}

\def\preuve{\begin{proof}}
\def\eop{\end{proof}}
\def\cqfd{\skip10=\parfillskip\parfillskip=0pt
\enspace\hfill\symbolecqfd\par\parfillskip=\skip10\par\medskip}
\def\symbolecqfd{\rlap{$\sqcap$}$\sqcup$}

\title{On the generalized membership problem in relatively hyperbolic groups}

\author{
  Olga Kharlampovich\thanks{Partially supported by the Simons Foundation, Award 422503.}\\
  \textit{\small Hunter College, CUNY}\\
  {\small\url{okharlampovich@gmail.com}}
  \and
  Pascal Weil\thanks{Partially supported by the DeLTA project (ANR-16-CE40-0007).}\\
 {\small \textit{Univ. Bordeaux, CNRS, Bordeaux INP, LaBRI, UMR5800, France}}\\
  {\small \textit{CNRS, ReLaX, UMI2000, Siruseri, India}}\\
  {\small \url{pascal.weil@labri.fr}}
}

\date{\today}

\begin{document}

\maketitle

\begin{abstract}
The aim of this note is to provide a proof of the decidability of the generalized membership problem for relatively quasi-convex subgroups of finitely presented relatively hyperbolic groups, under some reasonably mild conditions on the peripheral structure of these groups. These hypotheses are satisfied, in particular, by toral relatively hyperbolic groups.
\end{abstract}

\bigskip

The problem we consider here is the so-called \emph{generalized membership problem}, in a group $G$ generated by a finite set $A$: given a tuple $h_1, \ldots, h_k \in F(A)$ and letting $H$ be the subgroup they generate in $G$ (that is: $H$ is the subgroup of $G$ generated by the images of the $h_i$ in $G$), given an additional element $g\in G$ (also in the form of a word in $F(A)$), decide whether $g \in H$.

Stated as above, this problem is known to be undecidable without strong assumptions on the group $G$. Even in the relatively simple case of the direct product of two rank 2 free groups, $F_2 \times F_2$, there are finitely generated subgroups with undecidable membership problem (see Mihailova's subgroup \cite{1966:Mihailova}).

Our main result deals with the case where $G = \langle A \mid R\rangle$ is finitely presented and relatively hyperbolic with respect to a peripheral structure subject to additional conditions --- satisfied, in particular, by toral relatively hyperbolic groups, see Section~\ref{sec: main result}. Note that in these groups, and even in hyperbolic groups, there are finitely generated subgroups with undecidable membership problem \cite{1982:Rips}. We offer a \emph{partial algorithm} for the generalized membership problem in the following sense: an algorithm which may not stop on all instances but which will stop at least on those instances where $H$ is relatively quasi-convex, and which decides whether $g\in H$ when it stops.

We first survey some algorithmic results for groups, largely centered around this generalized membership problem, mainly focussing on those results that use graph-theoretic representations of subgroups, in particular the so-called Stallings graphs, as these are essential to our main result.

In the second section we use these results, and other results on the structure of relatively hyperbolic groups, to establish our main theorem.

%%%%%%%%%%%%%%
\section{Stallings graphs and algorithmic problems}

Stallings \cite{1983:Stallings} formalized a method (now known as \emph{Stallings foldings}) to associate with any finitely generated subgroup $H$ of a free group $F(A)$ an effectively computable discrete structure, called the \emph{Stallings graph} of $H$. This is a finite, oriented, labeled graph (the edges are labeled by elements of $A$) with a designated base vertex, in which the loops at the base vertex are labeled by reduced words representing the elements of $H$. Given a finite set of words in $F(A)$, we can compute the Stallings graph of the subgroup $H$ they generate (in time almost linear \cite{2006:Touikan}), compute the index, the rank and a basis for $H$, and solve the membership for $H$. In particular, this provides an elegant and computationally efficient solution of the generalized membership problem in $F(A)$: on input $(h_1,\ldots,h_n;g)$, one first computes the Stallings graph $\Gamma$ of the subgroup generated by the $h_i$, and one then verifies whether the reduced word $g$ can be read as a loop at the base vertex of $\Gamma$.

Given generators for another subgroup $K$, we can use the same tool of Stallings graphs to decide whether $H$ and $K$ are conjugates, compute their intersection and the finite collection (up to conjugacy) of intersections of their conjugates, and generally solve many other algorithmic problems, see \emph{e.g.} \cite{2002:KapovichMyasnikov,2007:MiasnikovVenturaWeil,2007:RoigVenturaWeil}. Most of these problems are solved very efficiently (in polynomial time) by this method, see \cite{2000:BirgetMargolisMeakin,2002:KapovichMyasnikov}.

Several authors introduced similar constructions to study finitely generated subgroups of non-free groups. More specifically, we are talking here of having an effectively constructible labeled graph canonically associated with a subgroup, solving at least the membership problem and allowing the computation of intersections. 

As mentioned in the introduction, one certainly needs to impose constraints on the group $G$. We also need to formulate assumptions on the subgroup $H\le G$. Indeed, even in good situations (\textit{e.g.} $G$ is automatic, or even hyperbolic), not every finitely generated subgroup has decidable membership problem \cite{1982:Rips}.

Pioneer work (published in 1996) came from two directions. Kapovich \cite{1996:Kapovich} used the Todd-Coxeter enumeration scheme to produce ever larger fragments of the Schreier (coset) graph $\Schreier(G,H)$, and showed that, if $G$ is geodesically automatic and $H$ is quasi-convex, one can decide when to stop this process and produce a Stallings-like graph to decide the membership problem for $H$. This yields a partial algorithm for the generalized membership problem, which halts exactly when the subgroup $H$ is quasi-convex. At the same time, Arzhantseva and Ol’shanskii \cite{1996:ArzhantsevaOlshanskii} studied a construction, starting with the Stallings graph of the subgroup $H_0 \le F(A)$ generated by $h_1,\ldots, h_n$, and enriching it by a combination of Stallings foldings and surgical additions of fragments of relators of $G$ (the so-called \emph{AO-moves}). For each integer $k \ge 1$, they identified a small cancellation property which 'almost always' holds (it is exponentially generic among the presentations with $r$ relators, $r$ fixed) under which every $k$-generated subgroup is quasi-convex, and such that their construction halts after a finite number of moves, and solves the membership problem for $H$.

In \cite{2005:McCammondWise}, McCammond and Wise also start from the Stallings graph of the subgroup $H_0 \le F(A)$, which they refine by so-called \emph{2-cell attachments}. They then use a geometric assumption (the \emph{perimeter reduction} hypothesis) on the complex representing the presentation $G = \langle A \mid R \rangle$, to show that their construction halts and produces a Stallings-like graph. They show that this geometric assumption holds in particular when $R$ consists of large powers, or under certain combinatorial conditions.

Kapovich \cite{1996:Kapovich} used his result to show that one can compute the quasi-convexity index of a (quasi-convex) subgroup. Arzhantseva's and Ol'shanskii's method was used to prove that, generically, a finitely presented group satisfies the Howson property \cite{1998:Arzhantseva}, see \cite{2000:Arzhantseva,2005:KapovichSchupp} for other applications. McCammond's and Wise's perimeter reduction hypothesis also leads to a number of algorithmic results, including the construction of Stallings-like graphs and the solution of the membership problem in large classes of presentations, many of which are locally quasi-convex (every finitely generated subgroup is quasi-convex) \cite{2005:McCammondWise}, see \cite{2003:Schupp} for other applications.

A common feature of these papers above is that they identify a method to `grow' a labeled graph, starting from the Stallings graph of a subgroup of the free group, and then exploit additional assumptions on both $G$ and $H$ to show that this growing process can be `terminated' at some point. 

It is natural, if we are going to rely on methods where words label paths in  graphs (which one can view as automata), to consider, as Kapovich \cite{1996:Kapovich} does, finitely presented groups $G = \langle A \mid R \rangle$ equipped with an automatic structure, providing in particular a rational language\footnote{A language is \emph{rational} (or \emph{regular}) if it is accepted by a finite state automaton.} of representatives for the elements of $G$, that is, a rational language $L$ over the alphabet $A \cup A\inv$, composed of reduced words, and such that $\mu(L) = G$ (where $\mu\colon F(A) \to G$ is the canonical onto morphism from the free group over $A$ onto $G$). It is also natural in this context to consider only so-called \emph{$L$-rational} subgroups $H$, that is, subgroups such that $L \cap \mu\inv(H)$ is a rational set as well. The notion of $L$-rationality, first considered by Gersten and Short \cite{1991:GerstenShort}, is equivalent to a geometric notion of $L$-quasi-convexity\footnote{Namely: there exists $\delta >0$ such that every $L$-representative of an element of $H$ stays within distance $\delta$ of $H$, in the Cayley graph of $G$.}. Classical quasi-convexity corresponds to the case where $L$ is the set of geodesic representatives of the elements of $G$. See Short \cite{1991:Short} for an example of the usage of automata-theoretic ideas to investigate quasi-convex subgroups.

An abundant literature considers the same set of problems for more specific classes of groups. Cai \emph{et al.} \cite{1994:CaiFuchsKozen} and later Gurevich and Schupp \cite{2007:GurevichSchupp} investigate the complexity of the generalized membership problem in the modular group. Schupp \cite{2003:Schupp} applies the results of \cite{2005:McCammondWise} to large classes of Coxeter groups, which turn out to be locally quasi-convex.  
Kapovich, Miasnikov, Weidmann \cite{2005:KapovichWeidmannMyasnikov} solve the membership
problem for subgroups of certain graphs of groups.
Markus-Epstein \cite{2007:Markus-Epstein} constructs a Stallings graph for the subgroups of amalgamated products of finite groups.
Silva, Soler-Escriva, Ventura \cite{2016:SilvaSoler-EscrivaVentura} do the same for subgroups of virtually free groups.
Here again, the groups considered are locally quasi-convex, and the authors rely on a folding process, much like in the free group case,
and a well-chosen set of representatives.
Finally, we mention Delgado and Ventura's work \cite{2013:DelgadoVentura}, where they develop a strong generalization of Stallings graphs to represent, and to compute with, subgroups of direct products of free and free abelian groups.

In \cite{2017:KharlampovichMiasnikovWeil}, Kharlampovich \emph{et al.} proposed a general approach to generalize a number of the situations listed above while keeping the spirit of the construction of Stallings graphs. If $G$ is an $A$-generated group and $L$ is a set of (possibly not unique) representatives for the elements of $G$, we define the \emph{Stallings graph of a subgroup $H$ with respect to $L$} to be the fragment $\Gamma_L(H)$ of the Schreier graph $\Schreier(G,H)$ spanned by the loops at vertex $H$ labeled by words of $L$ (that is: by the $L$-representatives of the elements of $H$). It is easily verified that this graph is finite if and only if $H$ is $L$-quasi-convex, or $L$-rational. We will use the following result in the next section.

\begin{theorem}[\cite{2017:KharlampovichMiasnikovWeil}]\label{thm: KhMW}
Let $G$ be an $A$-generated group, equipped with an automatic structure with language of representatives $L$. There exists a partial algorithm which, given $g,h_1,\ldots,h_k \in F(A)$, halts exactly if the subgroup $H$ generated by the $h_i$ is $L$-quasi-convex, and in that case outputs the Stallings graph of $H$ with respect to $L$.
\end{theorem}

Note that hyperbolic groups admit an automatic structure with language of representatives the set of geodesics. In particular, the partial algorithm in Theorem~\ref{thm: KhMW} computes a Stallings graph for the quasi-convex subgroups of hyperbolic groups.

Theorem~\ref{thm: KhMW} yields a uniform method to solve algorithmic problems for $L$-quasi-convex subgroups in automatic groups, including the generalized word problem and the computation of intersections. It also allows deciding conjugacy and almost malnormality, provided the automatic structure on $G$ satisfies a quantitative version of Hruska's and Wise's bounded packing property \cite{2009:HruskaWise} (this condition is satisfied by the geodesic automatic structure of hyperbolic groups), see \cite{2017:KharlampovichMiasnikovWeil}.

Theorem~\ref{thm: KhMW} was recently used by Kim in \cite{2019:Kim} where she, in particular, detects stability and Morseness in toral relatively hyperbolic groups

\begin{remark}(\cite{1996:Kapovich,2017:KharlampovichMiasnikovWeil})\label{rk: GMP from KhMW}
For the generalized membership problem in particular, the partial algorithm (halting exactly if $H$ is $L$-quasi-convex) consists in computing the Stallings graph $\Gamma$ as in Theorem~\ref{thm: KhMW}, using the automatic structure to compute an $L$-representative $w$ of $g$, and verifying whether $w$ labels a loop at the base vertex of $\Gamma$.
\end{remark}

%%%%%%%%%%%%%%%%%%
\section{The generalized membership problem for relatively hyperbolic groups}\label{sec: main result}

Let $G = \langle A \mid R\rangle$ be a finitely presented group and let $\calP$ be a finite collection of finitely generated subgroups of $G$, called the \emph{peripheral subgroups} of $G$. There are several definitions of $G$ being \emph{relatively hyperbolic with respect to the peripheral structure $\calP$}, due to Gromov \cite{1987:Gromov}, Farb \cite{1998:Farb}, Bowditch \cite{2012:Bowditch}, Dru\c tu and Sapir \cite{2005:DrutuSapir}, Osin \cite{2006:Osin}. These definitions turn out to be equivalent (see Bumagin \cite{2005:Bumagin}, Dahmani \cite{2003:Dahmani}, Hruska \cite[Theorem 5.1]{2010:Hruska}), we refer to the literature for details \cite{2006:Osin,2010:Hruska}.

If $H$ is a subgroup of $G$, there are also several definitions of relative quasi-convexity for $H$, in terms of natural geometries on $G$. Again, these are equivalent (Hruska \cite{2010:Hruska}) and we refer to the literature for precise definitions.

Properties of the \emph{parabolic subgroups} of $H$, that is, the subgroups that are contained in a conjugate of a peripheral subgroup $P\in \calP$, characterize certain subclasses of relatively quasi-convex subgroups, which will be useful in the sequel. We say that $H$ is \emph{peripherally finite} if every $H\cap P^x$ ($P\in \calP$, $x\in G$) is finite\footnote{These subgroups are called \emph{strongly quasi-convex} in \cite{2006:Osin}, and differ from the strongly quasi-convex subgroups of Tran \cite{2019:Tran}.}; more generally, we say that $H$ has \emph{peripherally finite index} if every infinite $H \cap P^x$ has finite index in $P^x$\footnote{These subgroups are called \emph{fully quasi-convex} in \cite{2010:ManningMarti-nez-Pedroza}}. Such subgroups are always finitely generated (Osin \cite[Thms 4.13 and 4.16]{2006:Osin} for the peripherally finite case, Kharlampovich \emph{et al.} \cite{2017:KharlampovichMiasnikovWeil} for the peripherally finite index case).

To go forward, we introduce the following assumptions on the peripheral structure $\calP$ of the relatively hyperbolic group $G$.

\paragraph{Assumptions (\textsf{Hyp})}
\begin{itemize}
\item[(H1)] Each group $P \in \calP$ satisfies the following: we are given a geodesically bi-automatic structure for $P$, on an alphabet $X_P$ and with language of representatives $L_P$, and we can compute a geodesically bi-automatic structure on every finite generating set of $P$ (given as a subset of $F(X_P)$).

\item[(H2)] The groups in $\calP$ are \emph{slender} (a.k.a. \emph{noetherian}: every one of their subgroups is finitely generated) and LERF.

\item[(H3)] For each $P\in\calP$, the set of tuples of words in $L_P$ that generate a finite index subgroup of $P$ is recursively enumerable.

\item[(H4)] We can solve the generalized membership problem in each $P\in \calP$.
\end{itemize}

\begin{remark}
Hruska showed that every relatively quasi-convex subgroup of $G$ is finitely generated, if and only if every group in $\calP$ is slender \cite[Cor. 9.2]{2010:Hruska}, so (H2) is a reasonable hypothesis to make in this algorithmic context.
\end{remark}

\begin{remark}
(\textsf{Hyp}) is satisfied in particular if the peripheral structure $\calP$ consists of finitely generated abelian groups, and notably, if $G$ is \emph{toral relatively hyperbolic} (that is: $G$ is torsion free and $\calP$ consists of non-cyclic free abelian groups).
\end{remark}

We can now state the central result of this note.

\begin{theorem}\label{thm: generalized membership}
Let $G = \langle A \mid R\rangle$ be a finitely presented group, relatively hyperbolic with respect to the peripheral structure $\calP$, and satisfying (\textsf{Hyp}). There is a partial algorithm which, given $g,h_1,\ldots,h_k \in F(A)$,
\begin{itemize}
\item halts at least if $g\in H$ or if the subgroup $H$ of $G$ generated by the $h_i$ is relatively quasi-convex and $g \not\in H$;

\item when it halts, decides whether $g\in H$.
\end{itemize}
\end{theorem}

The algorithm in Theorem~\ref{thm: generalized membership} is ``impractical'' in the following sense: there is no function bounding the time required for the algorithm to stop (if it will stop). It consists in two semi-algorithms, meant to be run concurrently, until one of them halts: one trying to witness the fact that $g\in H$ and the other trying to witness the opposite fact. 

The rest of this paper consists in the description of these semi-algorithms.

\paragraph{Semi-algorithm to verify that $g\in H$.}
It is a classical result that, given the presentation $\langle A\mid R\rangle$ for $G$ and given a word $g\in F(A)$, there is a partial algorithm which halts exactly if $g = 1$ in $G$. Indeed, $g = 1$ in $G$ if and only if a sequence of $R$-rewritings of $g$ eventually leads to the empty word. A systematic exploration of the $R$-rewritings of $g$ will eventually uncover this sequence if $g = 1$ in $G$.

This semi-algorithm is naturally extended to the problem at hand (does $g$ belong to $H$?) as follows. One starts with the Stallings graph $\Gamma$ of the subgroup of $F(A)$ generated by the $h_i$ (see \cite{1983:Stallings}), and iteratively:
\begin{itemize}
\item modify $\Gamma$ by gluing at every vertex a loop labeled by $r$ for every relator $r\in R$;

\item fold $\Gamma$ (this is the central step of the construction of Stallings graphs: it consists in identifying vertices $p$ and $q$ each time that there are edges labeled by a letter $a\in A$ from some vertex $s$ to both $p$ and $q$, or edges labeled by a letter $a$ from both $p$ and $q$ to some vertex $s$);

\item check whether $g$ labels a loop at the base vertex of $\Gamma$. If that is the case, then $g\in H$ and we stop the algorithm. If not, repeat.
\end{itemize}
A detailed discussion of this semi-algorithm can be found in \cite[Section 4.1]{2017:KharlampovichMiasnikovWeil}.

%%%%%%%%%%%%%%%%%%%%%%%%%
\paragraph{Semi-algorithm to verify that $g\not\in H$}
We call a subgroup of the form $H \cap P^x$ ($P\in\calP$, $x\in G$) which is infinite, a \emph{maximal infinite parabolic subgroup} of $H$. 
Our semi-algorithm relies on the following results.

\medskip

\noindent\textbf{[H]}\enspace Hruska shows \cite[Theorem 9.1]{2010:Hruska} that, if $H$ is relatively quasi-convex, then there exists a finite collection of maximal infinite parabolic subgroups $\{K_i\}_{1\le i \le \ell}$ such that every infinite maximal parabolic subgroup of $H$ is conjugated in $H$ to one of the $K_i$.

\medskip

\noindent\textbf{[MMP]}\enspace Manning and Mart\'\i nez-Pedroza show the following, under Hypothesis (H2) \cite[Theorem 1.7]{2010:ManningMarti-nez-Pedroza}. Suppose that $H \le G$ is relatively quasi-convex, $\{K_i\}_{1\le i \le \ell}$ is a collection of subgroups as in \textbf{[H]}, say with $K_i = H \cap P_i^{x_i}$ ($1\le i\le \ell$, $P_i\in \calP$, $x_i\in G$) and $g\not\in H$. Then there exist subgroups $R_i \le P_i^{x_i}$ such that $R_i$ has finite index in $P_i^{x_i}$, $K_i \le R_i$ and, if $K$ is generated by $H$ and the $R_i$, then $g\not\in K$ and $K$ has peripherally finite index. Note that \cite[Theorem 1.7]{2010:ManningMarti-nez-Pedroza} is a little more concise than this statement, which is extracted from the proof in that paper (\cite[p. 319]{2010:ManningMarti-nez-Pedroza}).

\medskip

\noindent\textbf{[AC]}\enspace Antolin and Ciobanu \cite[Cor. 1.9, Lemma 5.3, Thm 7.5]{2016:AntolinCiobanu} show that, under Hypothesis (H1), one can compute an automatic structure for $G$, with alphabet $X$ containing $A$ and the $X_P$ ($P\in\calP$), whose language $L$ of representatives consists only of geodesics (on alphabet $X$) and contains the $L_P$ ($P\in\calP$), and satisfying additional properties.

\medskip

\noindent\textbf{[KhMW]}\enspace Kharlampovich \emph{et al.} \cite[Sec. 7]{2017:KharlampovichMiasnikovWeil} build on \textbf{[AC]} to show that, if $H\le G$ is relatively quasi-convex (with respect to alphabet $A$) and has peripherally finite index, then it is $L$-quasi-convex with respect to alphabet $X$ \cite[Thm 7.5]{2017:KharlampovichMiasnikovWeil}. The proof of that theorem uses Hypothesis (H4). As explained in Remark~\ref{rk: GMP from KhMW}, this yields a solution of the membership problem in $H$.

\medskip

We can now give our semi-algorithm. For clarity, we give it as a non-deterministic partial algorithm. Such a non-determin\-istic algorithm can be turned into a deterministic one by standard methods (see, \textit{e.g.}, \cite[Thm 3.16]{2006:Sipser}).

\medskip

(1)\enspace We first apply \textbf{[AC]} to compute an automatic structure for $G$ on generator set $X$ (using Hypothesis (H1)). Then we compute a finite presentation of $G$ on $X$, say $\langle X \mid R_X\rangle$. For instance, $R_X$ consists of $R$, the relators $xu_x\inv$, where $x\in X\setminus A$ and $u_x$ is a fixed element of $F(A)$ such that $x = u_x$ in $G$, and all the cyclic permutations of these relators and their inverses.

The words $u_x$ can be computed as follows. Since the automatic structure for $G$ allows us to solve the word problem, one systematically checks whether $xu\inv$ is trivial, when $u$ runs through $F(A)$. As $G$ is  $A$-generated, some $u \in F(A)$ is equal to $x$ in $G$.

\medskip

(2)\enspace Choose non-deterministically a tuple $\vec x = (x_1,\cdots, x_\ell)$ of elements of $F(A)$; for each $1\le i \le \ell$, choose non-deterministically an element $P_i\in \calP$ and a tuple $\vec g_i$ of elements of $F(X_{P_i})$ generating a finite index subgroup of $P_i$ (this is possible under Hypothesis (H3)).

\medskip

(3)\enspace For this choice of $\vec x$ and the $\vec g_i$ ($1\le i\le \ell$), let $H_1 = \langle H \cup \bigcup_{i=1}^\ell{\vec g_i\/}^{x_i} \rangle$. Run the partial algorithm \textbf{[KhMW]} to decide whether $g\in H_1$ (using Hypothesis (H4)).

\medskip

Result \textbf{[MMP]} (which assumes Hypothesis (H2)), shows that, if $g\not\in H$ and $H$ is relatively quasi-convex, then for an appropriate choice of $\vec x$ and the $\vec g_i$, $H_1$ is relatively quasi-convex and has peripherally finite index, and $g\not\in H_1$. As $H_1$ has peripherally finite index, the partial  algorithm in Step (3) will halt and certify that $g\not\in H_1$, and hence that $g\not\in H$ since $H \le H_1$.

Summarizing: if $g\not\in H$ and $H$ is relatively quasi-convex, then one of the non-deterministic choices in Step (2) will be such that the partial algorithm halts and states that $g\not\in H$.  This completes the proof of Theorem~\ref{thm: generalized membership}.

%%%%%%%%%%%%%%%%%%%
%%%%%%%%%%%%%%%%%%%
%\def\cprime{$'$}
\bibliographystyle{abbrv}
\bibliography{pwbiblio}

%%%%%%%%%%%%%%%%%%%%%%%%
%%%%%%%%%%%%%%%%%%%%%%%%
%%%%%%%%%%%%%%%%%%%%%%%%
\end{document}